\documentclass[11pt]{amsart}
\usepackage{mathrsfs}
\usepackage{amssymb}

\usepackage{titletoc}
\usepackage{amscd}
\usepackage{amsmath, amssymb}
\usepackage{amsfonts}
\usepackage[colorlinks,linkcolor=blue,citecolor=blue, pdfstartview=FitH]{hyperref}

\numberwithin{equation}{section}

\theoremstyle{plain}
\newtheorem{thm}{Theorem}[section]
\newtheorem{theorem}[thm]{Theorem}
\newtheorem{lemma}[thm]{Lemma}
\newtheorem{corollary}[thm]{Corollary}
\newtheorem{proposition}[thm]{Proposition}

\theoremstyle{definition}

\newtheorem{question}[thm]{Question}

\newtheorem{remark}[thm]{Remark}

\newtheorem{definition}[thm]{Definition}

\newtheorem{example}[thm]{Example}

\newtheorem{defn-thm}[thm]{Definition-Theorem}

\newcommand{\sO}{{\mathcal O}}

\newcommand{\C}{{\mathbb C}}

\renewcommand{\H}{{\mathbb H}}

\newcommand{\R}{{\mathbb R}}

\renewcommand{\S}{{\mathbb S}}

\newcommand{\Z}{{\mathbb Z}}

\newcommand{\qtq}[1]{\quad\mbox{#1}\quad}
\newcommand{\bp}{\bar{\partial}}
\newcommand{\Om}{\Omega}

\newcommand{\ds}{\oplus}

\newcommand{\ts}{\otimes}

\newcommand{\btheorem}{\begin{theorem}}
\newcommand{\etheorem}{\end{theorem}}
\newcommand{\bproposition}{\begin{proposition}}
\newcommand{\eproposition}{\end{proposition}}
\newcommand{\bdefinition}{\begin{definition}}
\newcommand{\edefinition}{\end{definition}}
\newcommand{\bcorollary}{\begin{corollary}}
\newcommand{\ecorollary}{\end{corollary}}
\newcommand{\bproof}{\begin{proof}}
\newcommand{\eproof}{\end{proof}}
\newcommand{\bremark}{\begin{remark}}
\newcommand{\eremark}{\end{remark}}
\newcommand{\eexample}{\end{example}}
\newcommand{\bexample}{\begin{example}}
\newcommand{\la}{\langle}
\newcommand{\elemma}{\end{lemma}}
\newcommand{\blemma}{\begin{lemma}}
\newcommand{\ra}{\rangle}
\newcommand{\sq}{\sqrt{-1}}

\newcommand{\p}{\partial}

\renewcommand{\bar}{\overline}

\renewcommand{\phi}{\varphi}

\newcommand{\ee}{\end{eqnarray*}}
\newcommand{\be}{\begin{eqnarray*}}

\newcommand{\beq}{\begin{equation}}
\newcommand{\eeq}{\end{equation}}

\newcommand{\bd}{\begin{enumerate}}
\newcommand{\ed}{\end{enumerate}}

\renewcommand{\hat}{\widehat}

\renewcommand{\>}{\rightarrow}

\usepackage{fancyhdr}
\begin{document}
\title{Minimal complex surfaces with Levi-Civita Ricci-flat metrics}

 \makeatletter
\let\uppercasenonmath\@gobble
\let\MakeUppercase\relax
\let\scshape\relax
\makeatother
\author{Kefeng Liu}
\address{Department of Mathematics,
Capital Normal University, Beijing, 100048, China}
\address{Department of Mathematics, University of California at Los Angeles, California 90095}

\email{liu@math.ucla.edu}

\author{Xiaokui Yang}
\date{}
\address{{Address of Xiaokui Yang: Morningside Center of Mathematics, Institute of
        Mathematics, Hua Loo-Keng Key Laboratory of Mathematics,
        Academy of Mathematics and Systems Science,
        Chinese Academy of Sciences, Beijing, 100190, China.}}
\email{\href{mailto:xkyang@amss.ac.cn}{{xkyang@amss.ac.cn}}}

\maketitle \vskip -0.5cm {\hfill\small{\emph{In memory of Professor
Lu Qi-Keng\hfill}}

\begin{abstract}  This is a continuation of our previous paper \cite{LY14}.
In \cite{LY14}, we introduced the first Aeppli-Chern class on
compact complex manifolds, and proved that the  $(1,1)$ curvature
form of the Levi-Civita connection represents the first Aeppli-Chern
class which is a natural link between Riemannian geometry and
complex geometry. In this paper, we study the geometry of compact
complex manifolds with Levi-Civita Ricci-flat metrics and
 classify  minimal complex surfaces with
Levi-Civita Ricci-flat metrics. More precisely, we show that minimal
complex surfaces admitting Levi-Civita Ricci-flat metrics are
K\"ahler Calabi-Yau surfaces and Hopf surfaces.
\end{abstract}

\small{\setcounter{tocdepth}{1} \tableofcontents}

\section{Introduction}

 In this paper, we study the relationship between Riemannian
 manifolds and complex manifolds by using various metric connections and their curvature tensors.

 Let $(X,h)$ be a Hermitian manifold and $g$ be the background Riemannian metric.
 It is well-known that, when $(X,h)$ is not K\"ahler, the relation
 between the Riemannian geometry $(X,g)$ and the complex geometry
 $(X,h)$ is extremely complicated.  Indeed,
 on the Hermitian holomorphic tangent bundle $(T^{1,0}X,h)$, there are two
 typical metric compatible connections:
\bd \item the Chern connection $\nabla$, i.e. the unique connection
$\nabla$ compatible with the Hermitian metric and also the complex
structure $\bp$;

\item the Levi-Civita connection $\nabla^{\text{LC}}$, i.e. the
restriction of the complexified Levi-Civita connection on $T_\C X$
to the holomorphic tangent bundle $ T^{1,0}X$. \ed

\noindent From the definition,  it is quite obvious  that the
Levi-Civita connection $\nabla^{\text{LC}}$ is a representative of
the Riemannian geometry of  $(X,g)$. It is also well-known that when
$(X,h)$ is not K\"ahler, $\nabla$ and $\nabla^{\text{LC}}$ are not
the same. The complex geometry of the Chern connection is
extensively investigated in the literatures by using various methods
(e.g. \cite{CRS,F3, FY, FLY, L, LY12, LY14, ST2, ST3, STW, T,
TW,TW2, Yang16, Yang17A, Yang17, Yau2}). However, the complex
geometry of the Levi-Civita connection is not well understood
although it has rich Riemannian geometry
 structures..

 In \cite{LY14}, we introduced the first Aeppli-Chern classes for
 holomorphic line bundles. Let $L\to X$ be a holomorphic line bundle over $X$. The
first Aeppli-Chern class is defined as \beq
c^{\text{AC}}_1(L)=\left[-\sq\p\bp\log h\right]_{\text{A}} \in
H^{1,1}_{\text{A}}(X)\eeq where $h$ is an arbitrary smooth Hermitian
metric on $L$ and the  Aeppli cohomology is
  $$  H^{p,q}_{\text{A}}(X):=\frac{\text{Ker} \p\bp \cap \Om^{p,q}(X)}{\text{Im} \p \cap \Om^{p,q}(X)+ \text{Im}\bp \cap
  \Om^{p,q}(X)}.$$
 For a complex manifold $X$,
$c_1^{\mathrm{AC}}(X)$ is defined to be
$c_1^{\mathrm{AC}}(K^{-1}_X)$ where $K_X^{-1}$ is the anti-canonical
line bundle of $X$.  Note that, for a Hermitian line bundle $(L,h)$,
the classes $c_1(L)$ and $c_1^{\textrm{AC}}(L)$ have the same
$(1,1)$-form representative $\Theta^h=-\sq\p\bp\log h$ (in different
classes). It is well-known that on a Hermitian manifold
$(X,\omega)$, the first Chern-Ricci curvature
$\text{Ric}(\omega)=-\sq\p\bp\log\det(\omega)$ represents the first
Chern class $c_1(X)$. As an analog, we proved in
\cite[Theorem1.1]{LY14} that the first Levi-Civita Ricci curvature
$\mathfrak{Ric}(\omega)$ represents the first Aeppli-Chern class
$c_1^{\text{AC}}(X)$. It is obvious that $ c_1(X)=0 $ implies
$c_1^{\mathrm{AC}}(X)=0.$ Hence, it is very natural to study
non-K\"ahler Calabi-Yau manifolds by using the first Aeppli-Chern
class $c_1^{\text{AC}}(X)$ and the first Levi-Civita Ricci curvature
$\mathfrak{Ric}(\omega)$.  By the celebrated Calabi-Yau theorem
(\cite{Yau78}),   a compact K\"ahler manifold has $c_1(X)=0$ if and
only if it has a K\"ahler metric with Ricci-flat metric, i.e.
$\text{Ric}(\omega)=0$. It is easy to see that if
$\mathfrak{Ric}(\omega)=0$, then $c_1^{\text{AC}}(X)=0$. There is a
natural question analogous to the Calabi conjecture:
\begin{question}\label{question} On a compact complex manifold $X$, if   $c_1^{\text{AC}}(X)=0$, does there exist a smooth Levi-Civita
Ricci-flat Hermitian metric $\omega$, i.e.
$\mathfrak{Ric}(\omega)=0$?\end{question} \noindent As we have shown
in \cite[Theorem~1.2]{LY14}, $\mathfrak{Ric}(\omega)=0$ is
equivalent to \beq
\text{Ric}(\omega)=\frac{1}{2}(\p\p^*\omega+\bp\bp^*\omega).\label{MA}\eeq
The equation (\ref{MA}) is not the Monge-Amp\`ere type equation
since there are also non-elliptic second order derivatives on the
right hand side. As it is well-known,  it is particularly
challenging
 to solve such equations. Instead of solving the
equation (\ref{MA}) directly, we use several observations in
\cite{Yang17} to study the geometry of the equation (\ref{MA}) and
obtain necessary conditions to solve (\ref{MA}). See Corollary
\ref{top}, Theorem \ref{Kodaira} and Theorem \ref{ddbar} in Section
$5$ for more details. By using these necessary conditions and
Kodaira-Enriques' classification (e.g. \cite[p.~244]{BHPV}) of
minimal complex surfaces, we obtain:

\btheorem\label{cl} Let $X$ be a minimal complex surface. Suppose
$X$ admits a Levi-Civita Ricci-flat Hermitian metric $\omega$. Then
$X$ lies in one of the following \bd
\item  Enriques surfaces;\item
bi-elliptic surfaces;\item  K$3$ surfaces;
\item  $2$-tori;
\item  Hopf surfaces.
\ed \etheorem

\noindent It worths to point out that  we  do \emph{not} show every
Hopf surface  admitting a Levi-Civita Ricci-flat metric. We only
construct such metrics on  diagonal Hopf surfaces (see Theorem
\ref{hopf2}). We conjecture that all Hopf surfaces can support
Levi-Civita Ricci-flat metrics. On the other hand, every Hopf
surface $X$ has $c_1(X)=0\in H^2(X,\R)$. However, it is easy to show
that $X$ can not support Chern-Ricci flat Hermitian metrics, i.e.
Hermitian metrics $\omega$ with $\text{Ric}(\omega)=0$.\\

\noindent As an application of Theorem \ref{cl}, we obtain the
following example which  indicates that we need extra constraints to
solve  Question \ref{question} in general:

\bcorollary\label{example} Let $X$ be a Kodaira surface or an Inoue
surface. Then
$$c_1(X)=c_1^{\mathrm{BC}}(X)=c_1^{\mathrm{AC}}(X)=0.$$ However, $X$ does
not admit a Levi-Civita Ricci-flat Hermitian metric. \ecorollary

\noindent\textbf{Acknowledgement.} The second author would like to
thank  Valentino Tosatti for many useful comments and suggestions.
 This work
was partially supported   by China's Recruitment
 Program of Global Experts, National Center for Mathematics and Interdisciplinary Sciences,
 Chinese Academy of Sciences.

\section{Preliminaries}

\subsection{Chern connection on complex manifolds}
 Let $(X,\omega_g)$ be a compact Hermitian manifold. There exists a unique connection $\nabla$ on  the holomorphic tangent bundle $T^{1,0}X$
  which is compatible with the Hermitian metric and also the complex structure. This connection $\nabla$ is called the Chern connection.  The Chern
  connection $\nabla$
 on $(T^{1,0}X,\omega_g)$ has  curvature components \beq R_{i\bar
j k\bar \ell}=-\frac{\p^2g_{k\bar \ell}}{\p z^i\p\bar z^j}+g^{p\bar
q}\frac{\p g_{k\bar q}}{\p z^i}\frac{\p g_{p\bar \ell}}{\p\bar
z^j}.\eeq  The (first) Chern-Ricci form $\text{Ric}(\omega_g)$ of
$(X,\omega_g)$ has components
$$R_{i\bar j}=g^{k\bar \ell}R_{i\bar jk\bar \ell}=-\frac{\p^2\log\det(g)}{\p z^i\p\bar z^j}$$
which also represents the first Chern class $c_1(X)$ of the complex
manifold $X$. The Chern scalar curvature $s_g$ of $(X,\omega_g)$ is
given by \beq s_g=\text{tr}_{\omega_g}\text{Ric}(\omega_g)=g^{i\bar
j} R_{i\bar j}. \eeq The total scalar curvature of $\omega_g$ is
\beq \int_X s_g \omega_g^n=n \int\text{Ric}(\omega_g)\wedge
\omega_g^{n-1},\eeq where $n$ is the complex dimension of $X$.

\subsection{Bott-Chern classes and Aeppli classes}
The Bott-Chern
  cohomology and the Aeppli cohomology on a compact complex manifold
  $X$ are given by
  $$ H^{p,q}_{\mathrm{BC}}(X):= \frac{\text{Ker} d \cap \Om^{p,q}(X)}{\text{Im} \p\bp \cap \Om^{p,q}(X)}\qtq{and} H^{p,q}_{\mathrm{A}}(X):=\frac{\text{Ker} \p\bp \cap \Om^{p,q}(X)}{\text{Im} \p \cap \Om^{p,q}(X)+ \text{Im}\bp \cap
  \Om^{p,q}(X)}.$$

\noindent Let $\mathrm{Pic}(X)$ be the set of holomorphic line
bundles over $X$. As similar as the first Chern class map
$c_1:\mathrm{Pic}(X)\>H^{1,1}_{\bp}(X)$, there is a \emph{first
Aeppli-Chern class} map \beq c_1^{\mathrm{AC}}:\mathrm{Pic}(X)\>
H^{1,1}_{\mathrm{A}}(X).\eeq Given any holomorphic line bundle $L\to
X$ and any Hermitian metric $h$ on $L$, its curvature form
$\Theta_h$ is locally given by $-\sq\p\bp\log h$.  We define
$c_1^{\mathrm{AC}}(L)$ to be the class of $\Theta_h$ in
$H^{1,1}_{\mathrm{A}}(X)$.  For a complex manifold $X$,
$c_1^{\mathrm{AC}}(X)$ is defined to be
$c_1^{\mathrm{AC}}(K^{-1}_X)$ where $K_X^{-1}$ is the anti-canonical
line bundle. The first Bott-Chern class $c_1^{\mathrm{BC}}(X)$ can
be defined similary.

\subsection{Special manifolds} Let $X$ be a compact complex
manifold.

 \bd \item A Hermitian metric $\omega_g$ is called a Gauduchon metric if
$\p\bp\omega_g^{n-1}=0$. It is proved by Gauduchon (\cite{Ga2})
that, in the conformal class of each Hermitian metric, there exists
a unique Gauduchon metric (up to  scaling).

\item A Hermitian metric $\omega_g$ is called a balanced metric if
$d\omega_g^{n-1}=0$ or equivalently $d^*\omega_g=0$. On a compact
complex surface, a balanced metric is also K\"ahler, i.e.
$d\omega_g=0$. It is well-known many Hermitian manifolds can not
support balanced metrics, e.g. Hopf surface $\S^{3}\times \S^1$. It
is also obvious that balanced metrics are Gauduchon.

\item $X$ is called a Calabi-Yau manifold if $c_1(X)=0\in
H^2(X,\R)$.

\ed

\noindent It is obvious that,
$$c_1^{\mathrm{BC}}(X)=0 \Longrightarrow c_1(X)=0 \Longrightarrow c_1^{\mathrm{AC}}(X)=0,$$
and on compact K\"ahler manifolds or manifolds supporting the
$\p\bp$-lemma (\cite[Corollary~1.4]{LY14}), they are equivalent.

\vskip 2\baselineskip
\section{The Levi-Civita connection on the holomorphic tangent bundle}

Let's recall some elementary settings (e.g. \cite[Section~2]{LY14}).
Let $(M, g, \nabla)$ be a $2n$-dimensional Riemannian manifold with
the Levi-Civita connection $\nabla$. The tangent bundle of $M$ is
also denoted by $T_\R M$. The Riemannian curvature tensor of
$(M,g,\nabla)$ is $$
R(X,Y,Z,W)=g\left(\nabla_X\nabla_YZ-\nabla_Y\nabla_XZ-\nabla_{[X,Y]}Z,W\right)$$
for tangent vectors $X,Y,Z,W\in T_\R M$. Let $T_\C M=T_\R M\ts \C$
be the complexification. We can extend the metric $g$ and the
Levi-Civita connection $\nabla$ to $T_{\C}M$ in the $\C$-linear way.
Hence for any $a,b,c,d\in \C$ and $X,Y,Z,W\in T_\C M$, we have  $$
R(aX,bY,cZ, dW)=abcd\cdot R(X,Y,Z,W).$$

\noindent Let $(M,g,J)$ be an almost Hermitian manifold, i.e.,
$J:T_\R M\>T_\R M$ with $J^2=-1$, and
 for any $X,Y\in T_\R M$, $ g(JX,JY)=g(X,Y)$. The Nijenhuis tensor $N_J:\Gamma(M,T_\R M)\times \Gamma(M,T_\R
M)\>\Gamma(M,T_\R M)$ is defined as $$
N_J(X,Y)=[X,Y]+J[JX,Y]+J[X,JY]-[JX,JY].$$ The almost complex
structure $J$ is called \emph{integrable} if $N_J\equiv 0$ and then
we call $(M,g,J)$ a Hermitian manifold. We can also extend $J$ to
$T_\C M$ in the $\C$-linear way. Hence for any $X,Y\in T_\C M$, we
still have $ g(JX,JY)=g(X,Y).$
 By Newlander-Nirenberg's
theorem, there exists a real coordinate system $\{x^i,x^I\}$ such
that $z^i=x^i+\sq x^I$ are  local holomorphic coordinates on $M$.
Moreover, we have $T_\C M=T^{1,0}M\ds T^{0,1}M$ where
$$T^{1,0}M=\text{span}_\C\left\{\frac{\p}{\p z^1},\cdots,\frac{\p}{\p
z^n}\right\}\qtq{and}T^{0,1}M=\text{span}_\C\left\{\frac{\p}{\p \bar
z^1},\cdots,\frac{\p}{\p \bar z^n}\right\}.$$ Since $T^{1,0}M$ is a
subbundle of $T_{\C}M$, there is an induced connection
$\nabla^{\mathrm{LC}}$ on the holomorphic tangent bundle $T^{1,0}M$
given by \beq \nabla^{\text{LC}}=\pi\circ\nabla:
\Gamma(M,T^{1,0}M)\stackrel{\nabla}{\rightarrow}\Gamma(M, T_{\C}M\ts
T_{\C}M)\stackrel{\pi}{\rightarrow}\Gamma(M,T_{\C}M\ts T^{1,0}M).
\eeq Let $h=(h_{i\bar j})$ be the corresponding Hermitian metric on
$T^{1,0}M$ induced by $(M,g, J)$. It is obvious that
$\nabla^{\text{LC}}$ is a metric compatible connection on the
Hermitian holomorphic vector bundle $(T^{1,0}M, h)$, and we call
$\nabla^{\text{LC}}$ the \emph{Levi-Civita connection} on the
complex manifold $M$. It is obvious that, $\nabla^{\text{LC}}$ is
determined by the following relations \beq
\nabla^{\text{LC}}_{\frac{\p}{\p z^i}}\frac{\p}{\p
z^k}:=\Gamma_{ik}^p\frac{\p}{\p z^p} \qtq{and}
\nabla^{\mathrm{LC}}_{\frac{\p}{\p \bar z^j}}\frac{\p}{\p
z^k}:=\Gamma_{\bar jk}^p\frac{\p}{\p z^p} \eeq where \beq
\Gamma_{ij}^k=\frac{1}{2}h^{k\bar \ell}\left(\frac{\p h_{j\bar
\ell}}{\p z^i}+\frac{\p h_{i\bar \ell}}{\p z^j}\right), \qtq{and}
\Gamma_{\bar i j}^k=\frac{1}{2} h^{k\bar \ell}\left( \frac{\p
h_{j\bar\ell}}{\p\bar z^i}-\frac{\p h_{j\bar i}}{\p\bar
z^\ell}\right).\eeq The curvature  tensor $\mathfrak{R}\in
\Gamma(M,\Lambda^2 T_{\C}M\ts T^{*1,0}M\ts T^{1,0}M)$ of
$\nabla^{\mathrm{LC}}$ is given by $$ \mathfrak{R}(X,Y)s
=\nabla^{\mathrm{LC}}_{X}\nabla^{\mathrm{LC}}_Ys-\nabla^{\mathrm{LC}}_Y\nabla^{\mathrm{LC}}_Xs-\nabla^{\mathrm{LC}}_{[X,Y]}s$$
for any $X,Y\in T_{\C}M$ and $s\in T^{1,0}M$. A straightforward
computation shows that the curvature tensor $\mathfrak{R}$ has
$(1,1)$ components \beq \mathfrak{R}_{i\bar
jk}^{\ell}=-\left(\frac{\p \Gamma^{\ell}_{ik}}{\p \bar z^j}-\frac{\p
\Gamma^{\ell}_{\bar jk}}{\p z^i}+\Gamma_{ ik}^{s}\Gamma^{\ell}_{\bar
js}-\Gamma_{ \bar jk}^{s}\Gamma^{\ell}_{s i}\right).\label{curv}\eeq

 \bdefinition The \emph{(first) Levi-Civita Ricci curvature} $\mathfrak{Ric}(\omega_h)$  of
the Hermitian  vector bundle $\left(T^{1,0}M, \omega_h,
\nabla^{\mathrm{LC}}\right)$ is
 \beq \mathfrak{Ric}(\omega_h)=\sq\mathfrak
R^{(1)}_{i\bar j}dz^i\wedge d\bar z^j \qtq{with}
\mathfrak{R}^{(1)}_{i\bar j}=\mathfrak{R}_{i\bar j k}^{k}. \eeq  The
\emph{Levi-Civita scalar curvature} $ s_{\text{LC}}$ of
$\nabla^{\mathrm{LC}}$ on $T^{1,0}M$ is  \beq s_{\text{LC}}=h^{i\bar
j}h^{k\bar \ell}\mathfrak{R}_{i\bar j k\bar \ell}. \eeq
\edefinition

\vskip 2\baselineskip
\section{Geometry of the first Aeppli-Chern class}

Let's give a straightforward proof of  \cite[Theorem~1.2]{LY14}.

\btheorem\label{main-1} Let $(X,\omega)$ be a compact Hermitian
manifold. Then the first Levi-Civita Ricci form
$\mathfrak{Ric}(\omega)$ represents the first Aeppli-Chern class $
c_1^{\mathrm{AC}}(X)$ in $H^{1,1}_{\mathrm{A}}(X)$. Moreover, we
have the Ricci curvature relation \beq
\mathfrak{Ric}(\omega)=\emph{Ric}(\omega)-\frac{1}{2}(\p\p^*\omega+\bp\bp^*\omega),\label{key}
\eeq and  the scalar curvature relation \beq
s_{\mathrm{LC}}=s_{\mathrm{C}}-\la\p\p^*\omega,\omega\ra.\label{key2}\eeq
 \etheorem

\bproof It is easy to show that \beq \bp^*\omega=2\sq
\overline{\Gamma_{\bar{i}k}^k}dz^i \label{key11}\eeq and so
\beq-\frac{\p\p^*\omega+\bp\bp^*\omega}{2}=\sq\left(\frac{\p\Gamma_{\bar
j k}^k}{\p z^i}+\frac{\p\overline{\Gamma_{\bar{i}k}^k}}{\p\bar
z^j}\right)dz^i\wedge d\bar z^j.\eeq On the other hand, by formula
(\ref{curv}), we have

\beq  \mathfrak{R}_{i\bar j}=\mathfrak{R}_{i\bar jk}^k=-\frac{\p
\Gamma^{k}_{ik}}{\p \bar z^j}+\frac{\p \Gamma^{k}_{\bar jk}}{\p
z^i}.\eeq Moreover, we have \beq \left(-\frac{\p \Gamma^{k}_{ik}}{\p
\bar z^j}+\frac{\p \Gamma^{k}_{\bar jk}}{\p
z^i}\right)-\left(\frac{\p\Gamma_{\bar j k}^k}{\p
z^i}+\frac{\p\overline{\Gamma_{\bar{i}k}^k}}{\p\bar
z^j}\right)=-\frac{\p \Gamma^{k}_{ik}}{\p \bar
z^j}-\frac{\p\overline{\Gamma_{\bar{i}k}^k}}{\p\bar
z^j}=-\frac{\p^2\log \det(g)}{\p z^i\p\bar z^j}\eeq which
establishes formula (\ref{key}). It is easy to show
$\la\p\p^*\omega,\omega\ra=\la\bp\bp^*\omega,\omega\ra$.
 \eproof

\blemma\label{conformalchange} Let $(X,\omega)$ be a compact
Hermitian manifold with complex dimension $n$. Suppose $f\in
C^{\infty}(X,\R)$ and $\omega_f=e^f\omega$. Then we have \beq
 \bp^*_f\omega_f=\bp^*\omega+\sq(n-1)\p f\qtq{and} \bp\bp^*_f\omega_f= \bp\bp^*\omega-\sq(n-1)\p\bp f.\label{key12}\eeq
where $\bp^*, \bp^*_f$ are the adjoint operators with respect to the
metric $\omega$ and $\omega_f$ respectively. \elemma

\bdefinition The Kodaira dimension $\kappa(L)$ of a line bundle $L$
is defined to be
$$\kappa(L):=\limsup_{m\>+\infty} \frac{\log \dim_\C
H^0(X,L^{\ts m})}{\log m}$$ and the \emph{Kodaira dimension}
$\kappa(X)$ of $X$ is defined as $ \kappa(X):=\kappa(K_X)$ where the
logarithm of zero is defined to be $-\infty$.

\edefinition

\btheorem\label{scalarvanishing} Let $(X,\omega)$ be a compact
Hermitian manifold. Suppose the Levi-Civita scalar curvature
$s_{\mathrm{LC}}$ of $\omega$ is positive, then $K_X$ is not
pseudo-effective and $\kappa(X)=-\infty$. \etheorem

\bproof Let $\omega_f=e^{f}\omega$ be the Gauduchon metric
 in the conformal class of $\omega$. Then by Lemma
\ref{conformalchange} and Theorem \ref{main-1}, we have
\be\text{Ric}(\omega_f)-\frac{\p\p^*_f\omega_f+\bp\bp_f^*\omega_f}{2}&=&\text{Ric}(\omega)-\frac{\p\p^*\omega+\bp\bp^*\omega}{2}-\sq\p\bp
f\\
&=&\mathfrak{Ric}(\omega)-\sq\p\bp f. \ee

\noindent Moreover, we have
\begin{eqnarray} \int_X\text{Ric}(\omega_f)\wedge
\omega_f^{n-1}\nonumber&=&\int_X\left(\mathfrak{Ric}(\omega)-\sq\p\bp
f+\frac{\p\p^*_f\omega_f+\bp\bp^*\omega_f}{2}\right)\wedge
\omega_f^{n-1}\\
\nonumber&=&\int_X\mathfrak{Ric}(\omega)\wedge \omega_f^{n-1}+
\frac{1}{2}\left(\|\bp^*_f\omega_f\|^2_{\omega_f}+\|\p^*_f\omega_f\|^2_{\omega_f}\right)\\
&=&\frac{1}{n}\int_X e^{(n-1)f}\cdot s_{\mathrm{LC}}\cdot
\omega^{n}+\frac{1}{2}\left(\|\bp^*_f\omega_f\|^2_{\omega_f}+\|\p^*_f\omega_f\|^2_{\omega_f}\right).\label{Gauduchon}\end{eqnarray}
Suppose the Levi-Civita scalar curvature $s_{\mathrm{LC}}>0$, then
the total Chern scalar curvature of the Gauduchon metric $\omega_f$
is strictly positive, i.e.
$$\int_X s_{{f}}\cdot \omega_f^n=n\int_X\text{Ric}(\omega_f)\wedge \omega_f^{n-1}>0.$$
By \cite[Theorem~1.1]{Yang17} and \cite[Corollary~3.3]{Yang17}, we
know $K_X$ is not pseudo-effective and $\kappa(X)=-\infty$. \eproof

\bcorollary Let $X$ be a compact complex manifold with
$c_1^{\mathrm{BC}}(X)=0$, then there is no Hermitian metric with
positive Levi-Civita scalar curvature. \ecorollary

\bproof If $c_1^{\mathrm{BC}}(X)=0$, then by
\cite[Theorem~1.3]{STW}, there exists a smooth Gauduchon metric
$\omega_g$ with $\text{Ric}(\omega_g)=0$. Hence for any other
Gauduchon metric $\omega_G$, we have
$$\int_X \text{Ric}(\omega_G)\wedge \omega_G^{n-1}=0.$$
Suppose $\omega$ is a Hermitian metric with positive Levi-Civita
scalar curvature $s_{\mathrm{LC}}$ and $\omega_f$ is the Gauduchon
metric in the conformal class of $\omega$, then by
(\ref{Gauduchon}), we have

\begin{eqnarray*} \int_X\text{Ric}(\omega_f)\wedge
\omega_f^{n-1}=\frac{1}{n}\int_X e^{(n-1)f}\cdot
s_{\mathrm{LC}}\cdot
\omega^{n}+\frac{1}{2}\left(\|\bp^*_f\omega_f\|^2_{\omega_f}+\|\p^*_f\omega_f\|^2_{\omega_f}\right)>0,\end{eqnarray*}
which is a contradiction.
 \eproof

\vskip 2\baselineskip

\section{Compact complex manifolds with Levi-Civita Ricci-flat metrics}

 Let's recall that, a
Levi-Civita Ricci-flat metric is a Hermitian metric satisfying
$\mathfrak{Ric}(\omega)=0$, or equivalently, by formula
(\ref{key})\beq
\mathrm{Ric}(\omega)=\frac{\p\p^*\omega+\bp\bp^*\omega}{2}.\eeq

\noindent The first obstruction for the existence of Levi-Civita
Ricci-flat Hermitian metric is the top first Chern number:
\bcorollary\label{top} Suppose $c_1^{\mathrm{AC}}(X)=0$, then the
top intersection number $c_1^{n}(X)=0$. In particular, if $X$ has a
Levi-Civita Ricci-flat Hermitian metric $\omega$, then
$c_1^{n}(X)=0$. \ecorollary

\bproof By definition,  if $c_1^{\mathrm{AC}}(X)=0$, then
$$\text{Ric}(\omega)=\bp A+\p B$$
where $A$ is a $(1,0)$-form and $B$ is a $(0,1)$-form. Hence
$$c_1^n(X)=\int_X (\text{Ric}(\omega))^n=\int_X (\text{Ric}(\omega))^{n-1}\wedge (\bp A+\p B)=0$$
since $\text{Ric}(\omega)$ is both $\p$ and $\bp$-closed.
 \eproof

\btheorem\label{Kodaira} Let $X$ be a compact complex manifold.
Suppose $\omega$ is a  Levi-Civita Ricci-flat Hermitian metric. Then
either \bd \item $\kappa(X)=-\infty$; or\item  $\kappa(X)=0$ and
$(X,\omega)$ is conformally balanced with $K_X$ a holomorphic
torsion, i.e. $K_X^{\ts m}=\sO_X$ for some $m\in \Z^+$.\ed \etheorem

\bproof Let $\omega_f=e^{f}\omega$ be the Gauduchon metric in the
conformal class of $\omega$. Then by formula (\ref{Gauduchon}), the
total Chern scalar curvature of $\omega_f$ is \beq \int_X
s_{{f}}\cdot \omega_f^n=n\int_X\text{Ric}(\omega_f)\wedge
\omega_f^{n-1}=n\|\bp_f^*\omega_f\|_{\omega_f}^2\label{LCflat}\eeq
since the Levi-Civita scalar curvature $s_{\mathrm{LC}}=0$. Suppose
$\bp^*_f\omega_f\neq 0$,  then
$$\int_X s_{{f}}\cdot
\omega_f^n>0.$$ By \cite[Corollary~3.3]{Yang17}, we have
$\kappa(X)=-\infty$. On the other hand, if $\bp^*_f\omega_f=0$, i.e.
$(X,\omega)$ is conformally balanced. Then the total Chern scalar
curvature of the Gauduchon metric
$$\int_X s_{{f}}\cdot
\omega_f^n=0.$$ Then by \cite[Theorem~1.4]{Yang17}, we have
$\kappa(X)=-\infty$ or $\kappa(X)=0$, and when $\kappa(X)=0$, $K_X$
is a holomorphic torsion.
 \eproof

\bdefinition[\cite{DGMS75}] Let $X$ be a compact {complex} manifold.
$X$ is said to satisfy the $\p\bp$-lemma if the following statement
holds:
 if $\eta$ is $d$-exact, $\p$-closed and $\bp$-closed, it must be
 $\p\bp$-exact.
In particular, on such manifolds, for any pure-type form
$\phi\in\Om^{p,q}(X)$, if $\phi$ is $\bp$-closed and $\p$-exact,
then it is $\p\bp$-exact.\edefinition
 Let $\mu :
\hat X\>X$ be a modification between compact complex manifolds.  If
the $\p\bp$-lemma holds for $\hat X$, then the $\p\bp$-lemma also
holds for $X$. In particular, Moishezon manifolds  and manifolds in
Fujiki class $\mathscr C$  support the $\p\bp$-lemma. For more
details, we refer to \cite{DGMS75,AT13} and also the references
therein.

\btheorem\label{ddbar} Let $X$ be a compact complex manifold and
$\omega$ be a Levi-Civita Ricci-flat Hermitian metric. If $X$
supports the $\p\bp$-Lemma, then $(X,\omega)$ is conformally
balanced and $K_X$ is unitary flat. \etheorem

\bproof By formula (\ref{key}),  we have \beq
\text{Ric}(\omega)=\frac{\p\p^*\omega+\bp\bp^*\omega}{2}, \eeq since
$\mathfrak{Ric}(\omega)=0$. Note that $\text{Ric}(\omega)$ is
$\p$-closed and $\bp$-closed, and so we have
$$\bp\p\p^*\omega=0$$
 Moreover, if $X$ supports the
$\p\bp$-Lemma, then the $\bp$-closed and $\p$-exact $(1,1)$-form
$\p\p^*\omega$ is $\p\bp$-exact, i.e. there exists a smooth function
$\phi$ such that
$$\p\p^*\omega=\p\bp\phi.$$
Therefore, $\text{Ric}(\omega)=\sq\p\bp F$ where
$F=-\frac{\bar\phi-\phi}{2\sq}\in C^\infty(X,\R)$. It is obvious
that the Hermitian metric $e^{\frac{F}{n}}\omega$ is Chern
Ricci-flat, i.e. $K_X$ is unitary flat. Moreover, for any Gauduchon
metric $\omega_G$ on $X$, we have \beq \int_X
\text{Ric}(\omega_G)\wedge \omega_G^{n-1}=\int_X
\text{Ric}(\omega)\wedge \omega_G^{n-1}=\int_X\sq\p\bp F\wedge
\omega_G^{n-1}=0.\eeq Let $\omega_f=e^{f}\omega$ be the Gauduchon
metric  in the conformal class of $\omega$, then by (\ref{LCflat}),
we have
$$ \int_X\text{Ric}(\omega_f)\wedge
\omega_f^{n-1}=\|\bp_f^*\omega_f\|_{\omega_f}^2=0,$$ that is
$\bp^*_f\omega_f=0$.
 \eproof

\vskip 2\baselineskip
\section{Classification of minimal complex surfaces with Levi-Civita Ricci-flat metrics}

In this section, we investigate  the Levi-Civita Ricci-flat metrics
on minimal complex surfaces and prove Theorem \ref{cl}.

\btheorem\label{greater} Let $X$ be a minimal complex surface with
$\kappa(X)\geq 0$. Suppose $X$ admits a Levi-Civita Ricci-flat
Hermitian metric $\omega$. Then $X$ is K\"ahler surface of
Calabi-Yau type, i.e. $X$ is exactly one of the following \bd\item a
Enriques surface;\item a bi-elliptic surface;\item a K$3$ surface;
\item a torus.\ed \etheorem

\bproof By Theorem \ref{Kodaira}, we have $\kappa(X)\leq 0$. Hence,
we only need to consider minimal surfaces with $\kappa(X)=0$. By the
Kodaira-Enriques' classification of minimal surfaces (e.g. \cite[p.
244]{BHPV}), a minimal surface with $\kappa(X)=0$ has torsion
canonical line bundle $K_X$, i.e. $K_X^{\ts m}=\sO_X$ for some $m\in
\Z^+$. Hence, there exists a Hermitian metric $\omega_0$ with
$\text{Ric}(\omega_0)=0$. Let $\omega$ be the Hermitian metric with
$\mathfrak{Ric}(\omega)=0$ and $\omega_f=e^{f}\omega$ be the
Gauduchon metric  in the conformal class of $\omega$. Then \beq
\int_X \text{Ric}(\omega_f)\wedge \omega_f=\int_X
\text{Ric}(\omega_0)\wedge \omega_f=0.\eeq Hence, by (\ref{LCflat}),
we obtain $\bp_f^*\omega_f=0$. Since $\dim X=2$, we have
$d\omega_f=0$, i.e. $X$ is a K\"ahler surface. According to the
Kodaira-Enriques' classification, $X$ is either an Enriques surface,
a bi-elliptic surface, a K$3$ surface or a torus. All these surfaces
are K\"ahler surfaces of Calabi-Yau type, and all K\"ahler
Calabi-Yau metrics are Levi-Civita Ricci-flat. \eproof

\bremark The Hermitian metric $\omega$ with
$\mathfrak{Ric}(\omega)=0$ in Theorem \ref{greater}  is not
necessarily K\"ahler. Indeed, let $\omega_{\mathrm{CY}}$ be a
Calabi-Yau K\"ahler metric on $X$.
 Then for any non constant smooth function $f\in C^\infty(X,\R)$, we can construct a non-K\"ahler
 Levi-Civita Ricci-flat metric. By Yau's  theorem, there
 exists a K\"ahler metric $\omega_0$ such that
 $$\omega_0^2=e^{-f}\omega^2_{\mathrm{CY}}.$$
Let $\omega=e^f\omega_0$. Then $\omega$ is a non-K\"ahler metric
with Levi-Civita Ricci-flat curvature. Indeed, \be
\mathfrak{Ric}(\omega)&=&\text{Ric}(\omega)-\frac{\p\p^*\omega+\bp\bp^*\omega}{2}\\
&=&\text{Ric}(\omega_0)-2\sq\p\bp
f-\frac{\p\bp_0^*\omega_0+\bp\bp^*_0\omega_0}{2}+\sq\p\bp f\\
&=&\text{Ric}(\omega_{\mathrm{CY}})+\sq\p\bp f-2\sq\p\bp
f-\frac{\p\bp_0^*\omega_0+\bp\bp^*_0\omega_0}{2}+\sq\p\bp f \\
&=&0,\ee where we use Lemma \ref{conformalchange} in the second
identity.
 \eremark

\btheorem\label{infty} Let $X$ be a minimal complex surface with
$\kappa(X)=-\infty$. Suppose $X$ admits a Levi-Civita Ricci-flat
Hermitian metric $\omega$. Then $X$ is a Hopf surface. \etheorem

\bproof According to the Kodaira-Enriques' classification
\cite[p.~244]{BHPV}, $X$ is one of the following \bd \item a minimal
rational surface;
\item a minimal surface of class $\mathrm{VII}$;

\item a ruled surface of genus $g\geq 1$.

\ed Suppose $X$ is a K\"ahler surface, i.e. $X$ is a minimal
rational surface or a ruled surface of genus $g\geq 1$. Then the
$\p\bp$-Lemma holds on $X$. Hence by Theorem \ref{ddbar}, $K_X$ is
unitary flat and $\kappa(X)=0$ which a is contradiction. Hence $X$
is non-K\"ahler, i.e. $X$ is a minimal surface of class
$\mathrm{VII}$. By Corollary \ref{b2}, Theorem \ref{Inoue2} and
Theorem \ref{hopf2} in the next section, we complete the proof.
\eproof

\bcorollary Let $X$ be a Kodaira surface or an Inoue surface. Then
$$c_1^{\mathrm{BC}}(X)=c_1^{\mathrm{AC}}(X)=0.$$ However, $X$ does
not admit a Levi-Civita Ricci-flat Hermitian metric. \ecorollary

\bproof It is well-known (\cite[p.~244]{BHPV}) that a Kodaira
surface is a non-K\"ahler surface with torsion canonical line bundle
$K_X$. Hence $c_1^{\mathrm{BC}}(X)=c_1^{\mathrm{AC}}(X)=0$ and
$\kappa(X)=0$. By Theorem \ref{greater}, it has no Levi-Civita
Ricci-flat Hermitian metric. For an Inoue surface, one has
$b_2(X)=0$. Hence, $c_1(X)=c_1^{\mathrm{AC}}(X)=0$.\eproof

\vskip 2\baselineskip

\section{The Levi-Civita Ricci-flat metrics on minimal surfaces of class $\mathrm{VII}$}

A class $\mathrm{VII}$ surface is a minimal compact complex surface
with $b_1=1$ and $\kappa(X)=-\infty$. There are three classes of
them

\bd
\item Hopf surfaces: whose universal cover is $\C^2-\{0\}$, or
equivalently a class $\mathrm{VII}$ surface with $b_2=0$ and
contains a curve;

\item Inoue surfaces: a class $\mathrm{VII}$ surface has $b_2=0$ and
contains no curves;

\item all class $\mathrm{VII}$ surfaces with $b_2>0$.
\ed

\bcorollary\label{b2} On  $\mathrm{VII}$ surfaces with $b_2>0$,
there is no Levi-Civita Ricci-flat Hermitian metrics. \ecorollary

\bproof It is well-known that on  minimal $\mathrm{VII}$ surfaces we
have
$$c_1^2(X)=-b_2.$$
Corollary \ref{b2} follows from Corollary \ref{top}. \eproof

\subsection{Inoue surfaces} It is well-known (\cite{In}) that an
Inoue surface is a quotient of $\H\times \C$ by a properly
discontinuous group of affine transformations where $\H$ is the
upper half-plane. There are three types of Inoue surfaces:

\bd \item Inoue surfaces $S_M$. Let $M$ be a matrix in
$\mathrm{SL}_3(\Z)$ admitting one real eigenvalue $\alpha>1$ and two
complex conjugate eigenvalues $\beta\neq \bar\beta$. Let
$(a_1,a_2,a_3)$ be a real eigenvector of $M$ corresponding to
$\alpha$ and let $(b_1,b_2,b_3)$ be an eigenvector of $M$
corresponding to $\beta$. Then $X=S_M$ is the quotient of $\H\times
\C$ by the group of affine automorphisms generated by
$$g_0(w,z)=(\alpha w ,\beta z),$$
$$g_i(w,z)=(w+a_i,z+b_i), \ \ i=1,2,3.$$

\item Inoue surfaces $X=S^+_{N,p,q,r;t}$ are defined as the
quotient of $\H\times \C$ by the group of affine automorphisms
generated by
$$g_0(w, z)=(\alpha w, z+t),$$
$$g_i(w,z)=(w+a_i, z+b_iw+c_i),\ \ i=1,2$$
$$g_3(w,z)=\left(w,z+\frac{b_1a_2-b_2a_1}{r}\right),$$
where $(a_1,a_2)$ and $(b_1,b_2)$ are the eigenvectors of some
matrix $N\in \mathrm{SL}_2(\Z)$ admitting real eigenvalues
$\alpha>1$, $\alpha^{-1}$. Moreover $t\in \C$ and $p, q, r (r\neq
0)$ are integers, and $(c_1,c_2)$ depends on $(a_i,b_i), p, q, r$
(see \cite{In}).

\item Inoue surfaces $X=S^-_{N,p,q,r;t}$  have unramified double cover which are Inoue surfaces of
type $S^+_{N,p,q,r;t}$.

\ed

\btheorem\label{Inoue2} On Inoue surfaces, there is no Levi-Civita
Ricci-flat Hermitian metrics. \etheorem

\bproof Suppose $\omega$ is a Levi-Civita Ricci-flat Hermitian
metric on the Inoue surface $X$.  Let $\omega_f=e^{f}\omega$ be the
Gauduchon metric
 in the conformal class of $\omega$, then by formula (\ref{Gauduchon}), the
total Chern scalar curvature of $\omega_f$ is \beq \int_X
s_{{f}}\cdot \omega_f^2=2\int_X\text{Ric}(\omega_f)\wedge
\omega_f=2\|\bp_f^*\omega_f\|_{\omega_f}^2\geq 0.\label{Inoue1}\eeq
We shall show that on each Inoue surface, there exists a smooth
Gauduchon metric with non-positive but not identically zero first
Chern-Ricci curvature. Indeed, let $(w,z)\in \H\times \C$ be the
holomorphic coordinates, then by the precise definition of each
Inoue surface  (see also  \cite{DPS, FTWZ, Tel}), we know the form
$$\sigma=\frac{dw\wedge dz}{\text{Im}(w)}$$
descends to a smooth  nowhere vanishing $(2,0)$ form on $X$, i.e.
$\sigma\in \Gamma(X,K_X)$. Then it induces a smooth Hermitian metric
$h$ on $K_X$ given by $h(\sigma,\sigma)=1$.
 In the holomorphic frame $e=dw\wedge dz$ of $K_X$, we have
 $$h=h(e,e)=[\text{Im}(w)]^2.$$
It also induces a Hermitian metric $h^{-1}$ on $K_X^{-1}$, and the
curvature of $h^{-1}$ is
$$-\sq \p\bp\log h^{-1}=\sq\p\bp\log [\text{Im}(w)]^2=-\frac{\sq}{2}\frac{dw\wedge d\bar w}{[\text{Im}(w)]^2},$$
which also represents $c^{\mathrm{BC}}_1(X)$. By Theorem
\cite[Theorem~1.3]{STW}, there exists a Gauduchon metric $\omega_G$
with $$\text{Ric}(\omega_G)=-\frac{\sq}{2}\frac{dw\wedge d\bar
w}{[\text{Im}(w)]^2}\leq 0.$$ Hence, for any Gauduchon metric
$\omega$, one has
$$\int_X\text{Ric}(\omega)\wedge \omega=\int_X\text{Ric}(\omega_G)\wedge \omega<0$$
 which is a
contradiction to (\ref{Inoue1}). \eproof

\subsection{Hopf manifolds} Let's recall an example in \cite[Section~6]{LY12, LY14}.
Let $X=\S^{2n-1}\times \S^1$ be the standard $n$-dimensional ($n\geq
2$) Hopf manifold. It is diffeomorphic to $\C^n- \{0\}/G$ where $G$
is cyclic group generated by the transformation $z\rightarrow
\frac{1}{2}z$. It has  a naturally induced metric $\omega_0$ given
by
 \beq \omega_0=\sq
\frac{\delta_{i\bar j}}{|z|^2}dz^i\wedge  d\bar z^j. \eeq

\noindent We present a straightforward computation to show (c.f.
\cite[Theorem~6.2]{LY14}): \btheorem\label{hopf2} The perturbed
metric \beq \omega_g=\omega_0-\frac{1}{n}\cdot\sq\p\bp\log
|z|^2.\eeq  is Levi-Civita Ricci-flat, i.e. $
\mathfrak{Ric}(\omega_g)=0.$ \etheorem

 \bproof  If we write $ \omega_g=\sq  g_{i\bar j}dz^i\wedge
d\bar z^j$, then
$$ g_{i\bar j}=\frac{1}{|z|^2}\left(\frac{n-1}{n}\delta_{ij}+\frac{ \bar z^i z^j}{n|z|^2}\right), \qtq{and } g^{i\bar j}=|z|^2\left(\frac{n\delta_{i\bar
j}}{n-1}-\frac{ z^i\bar z^j}{(n-1)|z|^2}\right).$$ Let $\p^*$ and
$\bp^*$ be the adjoint operators with respect to $\omega_g$ and
$\Lambda$ is  the dual operator of $\omega_g\wedge\bullet$. A
straightforward  computation (\cite[Lemma~3.3]{LY14}) shows
$$\bp^*\omega_g=\sq\Lambda\p\omega_g.$$
Note  that $$\p\omega_g=\p\omega_0=-\frac{\sq \delta_{i\bar j}\bar
z^k}{|z|^4}dz^k\wedge dz^i\wedge d\bar z^j.$$ Hence, we obtain \be
\bp^*\omega_g=\sq\Lambda\p\omega_g&=&\sq g^{k\bar
q}\frac{\delta_{i\bar q}\bar z^k}{|z|^4}dz^i-\sq g^{i\bar
q}\frac{\delta_{i\bar q}\bar z^k}{|z|^4}dz^k\\
&=&\sq \frac{\sum_kg^{k\bar i}\bar z^k}{|z|^4}dz^i-\sq \frac{\sum_q
g^{q\bar q} \bar z^k}{|z|^4}dz^k.\ee On the other hand, \beq \sum_k
g^{k\bar i}\bar z^k=|z|^2\bar z^i \qtq{and} \sum_q g^{q\bar q}\bar
z^k=(n+1)|z|^2\bar z^k.\eeq Hence, we have
$$\bp^*\omega_g=\sq n \frac{\bar z^kdz^k}{|z|^2}=-n\sq \p\log |z|^2$$
and
$$\bp\bp^*\omega_g=n\sq\p\bp\log|z|^2.$$
Therefore
$$ \frac{\p\p^*\omega_g+\bp\bp^*\omega_g}{2}=
\sq n\p\bp\log |z|^2.$$  A direct computation shows $\det( g_{i\bar
j})=(1+\lambda)^{n-1}|z|^{-2n}$, and we have  $$
\text{Ric}(\omega_g)=-\sq\p\bp\log\det(g)=n\cdot \sq \p\bp\log
|z|^2.$$ By Theorem \ref{main-1},
$$\mathfrak{Ric}( \omega_g)=\text{Ric}( \omega_g)-\frac{\p\p^*\omega_g+\bp\bp^*\omega_g}{2}=0.$$

 \eproof

\bremark In this example, we construct a solution to the Levi-Civita
Ricci-flat equation on Hopf manifolds. It is natural to ask whether
there are more solutions. We expect there are theoretical approaches
on the existence of Levi-Civita Ricci-flat metrics on all Hopf
manifolds.
 \eremark


\begin{thebibliography}{99}

\bibitem{AT13} Angella, D.; Tomassini, A. On the $\p\bp$-lemma and Bott-Chern cohomology. Invent. Math. \textbf{192} (2013), no. \textbf{1}, 71--81.


\bibitem{BHPV} Barth, W.; Hulek, K.;  Peters, C.;  Van de Ven, A. \textit{Compact complex surfaces}.  Ergebnisse der Mathematik und ihrer Grenzgebiete. 3. Folge. A Series of Modern Surveys in Mathematics. Springer-Verlag, Berlin, 2004.

\bibitem{CRS} Chiose, I.; Rasdeaconu, R.; Suvaina, I. Balanced metrics on uniruled
manifolds.
 \href{http://arxiv.org/abs/1408.4769}{arXiv:1408.4769}. To appear
 in Comm. Anal. Geom.
\bibitem{DGMS75} Deligne, P.; Griffiths, P.; Morgan, J.; Sullivan, D.  Real homotopy theory of K\"ahler manifolds. Invent. Math. \textbf{29} (1975), no. 3, 245--274.




\bibitem{DPS}  Demailly, J.-P.; Peternell, T.; Schneider, M. Compact complex manifolds with numerically effective tangent bundles. J. Algebraic Geom. \textbf{3} (1994), no. 2, 295--345.


\bibitem{FTWZ} Fang, S.-W.;  Tosatti, V.;  Weinkove, B. and  Zheng, T.  Inoue surfaces and the Chern-Ricci
flow.  J. Funct. Anal. \textbf{271} (2016),  3162--3185.

 \bibitem{F3} Fu, J.-X.  On non-K\"ahler Calabi-Yau threefolds with balanced metrics. \emph{Proceedings of the International Congress of Mathematicians.} Volume II, 705--716, Hindustan Book Agency, New Delhi, 2010.
\bibitem{FY} Fu, J.- X.; Yau, S.-T. The theory of superstring with flux on non-K\"ahler manifolds and the complex Monge-Amp\'ere equation. J. Differential Geom. \textbf{78} (2008), no. 3, 369--428.
\bibitem{FLY}Fu, J.- X; Li, J.; Yau, S.-T. Constructing balanced metrics on some families of non-K\"ahler Calabi-Yau threefolds.  J. Differential Geom. \textbf{90} (2012), no. 1,81--129.

\bibitem{Ga2}  Gauduchon, P. { Fibr\'es hermitiens \`a endomorphisme de Ricci non-n\'egatif}, Bull. Soc. Math. France {\textbf {105}}  1977), 113--140.






\bibitem{In} Inoue, M. On surfaces of type $\mathrm{VII_0}$, Inv. Math., \textbf{24} (1974), 269--310.


\bibitem{L} Li, Yi. A priori estimates for Donaldson’s equation over compact Hermitian manifolds. Calc. Var. Partial Differential Equations, \textbf{50}(2014), no. 3-4, 867--882.

\bibitem{LY12} Liu, K.-F.; Yang, X.-K. Geometry of Hermitian manifolds.  Internat. J. Math. \textbf{23} (2012) 40pp.

\bibitem{LY14} Liu, K.-F.; Yang, X.-K. Ricci curvatures on Hermitian manifolds.{Trans. Amer. Math. Soc.}  \textbf{369} (2017),
5157--5196.




\bibitem{ST2} Streets, J.; Tian, G. {A parabolic flow of pluriclosed metrics}, Int. Math. Res. Not.  2010, no. 16, 3101--3133.
\bibitem{ST3} Streets, J.; Tian, G. {Regularity results for pluriclosed flow}, Geom. Topol. \textbf{17} (2013), no. 4, 2389--2429.


\bibitem{STW} Sz\'{e}kelyhidi, G.;  Tosatti, V.; Weinkove, B. Gauduchon metrics with prescribed volume
form. \href{http://arxiv.org/abs/1503.04491}{arXiv:1503.04491}.

\bibitem{Tel} Teleman, A. The pseudo-effective cone of a non-K\"ahlerian surface
and applications. Math. Ann. \textbf{335}(2006), 965--989.

\bibitem{T} Tosatti, V. Non-K\"ahler Calabi-Yau manifolds.  Contemp. Math. \textbf{644}(2015), 261--277.

\bibitem{TW} Tosatti, V.; Weinkove, B. On the evolution of a Hermitian metric by its Chern-Ricci
form. J. Differential Geom. \textbf{99} (2015), no.1, 125--163.

\bibitem{TW2} Tosatti, V.; Weinkove, B. The
Monge-Amp\`{e}re equation for $(n-1)$-plurisubharmonic functions on
a compact K\"ahler manifold. J. Amer. Math. Soc. 30 (2017), no.2,
311-346.



\bibitem{Yang16}Yang, X.-K. Hermitian manifolds with semi-positive
holomorphic sectional curvature.  {Math. Res. Lett.} \textbf{23}
(2016), no.3, 939--952.

\bibitem{Yang17A} Yang, X.-K.  The Chern-Ricci
flow and holomorphic bisectional curvature.
 Sci. China Math. \textbf{59} (2016), 2199-2204.

\bibitem{Yang17} Yang, X.-K. Scalar curvature on compact complex
manifolds. \href{https://arxiv.org/abs/1705.02672}{arXiv:1705.02672}

\bibitem{Yau2} Yau, S.-T. On the curvature of compact Hermitian manifolds. Invent. Math. \textbf{25} (1974), 213--239.
\bibitem{Yau78} Yau, S.-T.  { On the Ricci curvature of a compact K\"ahler manifold and the complex Monge-Amp\`ere equation, I},
Comm. Pure Appl. Math. {\textbf{ 31}} (1978),  339--411.

\end{thebibliography}
\end{document}